\documentclass[11pt]{article}
\usepackage{amssymb}
\usepackage{epsfig}
\usepackage{mathrsfs}

\def\eps{\varepsilon}
\def\ZZ{\mathbb Z}
\def\RR{\mathbb R}

\def\na{\nabla}
\def\pa{\partial}

\def\ds{\displaystyle}
\def\vphi{\varphi}

\def\De{{D_{\eps}}}
\def\ov{\overline v}
\def\ou{\overline u}

\def\ovae{{\overline v}_{\alpha}^\eps}
\def\vae{{ v}_{\alpha}^\eps}

\def\wae{\overline {w}_{\alpha}^\eps}

\def\om{\overline m}

\def\oell{\overline \ell}

\def\cor{\overline w^\eps}

\def\og{\overline \gamma}
\def\ug{\underline \gamma}

\def\J{\mathscr J}
\def\Ja{\mathscr J_\alpha}
\def\capa{{\rm cap}}

\def\A{D}

\newtheorem{theorem}{Theorem}[section]
\newtheorem{lemma}[theorem]{Lemma}
\newtheorem{proposition}[theorem]{Proposition}

\newtheorem{corollary}[theorem]{Corollary}

\def\qed{\hbox{${\vcenter{\vbox{                        
   \hrule height 0.4pt\hbox{\vrule width 0.4pt height 6pt
   \kern5pt\vrule width 0.4pt}\hrule height 0.4pt}}}$}}

\title{Random Homogenization of an Obstacle Problem}

\author{L. A. Caffarelli\footnote{Dept of Mathematics, University of Texas at Austin,  Austin, TX~78712, USA} \quad
and A. Mellet\footnote{Dept of Mathematics, University of British Columbia, Vancouver, BC V6T 1Z2, Canada}}

\date{\today}

\begin{document}

\maketitle

 \begin{abstract}
We study the homogenization of an obstacle problem in a perforated domain, when the holes have random shape and size. 
The main assumption concerns the capacity of  the holes which is assumed to be stationary ergodic.
 \end{abstract}

\section{Introduction}

Let $(\Omega,\mathcal F, \mathcal P)$ be a given probability space. For  every $\omega\in\Omega$ and every $\eps>0$, we consider a domain
$\De(\omega)$  obtained by perforating holes from a bounded domain $D$ of $\RR^{n}$.
We are interested in the asymptotic behavior as $\eps\rightarrow 0$ of the solution of the following obstacle problem:
$$ \min \left\{  \int_{D} \frac{1}{2} |\na u|^2 - f \, u\, dx\, ;\, u \geq 0   \mbox{ a.e. in } D \setminus \De\, , \; u \in  H^1_0 (D) \right\}
$$
for some $f\in L^2(D)$.
This is a well known homogenization problem and the asymptotic behavior of the solutions strongly depends on the size and the repartition of the holes 
$$T_{\eps}(\omega)= D\setminus\De.$$
This problem was first studied in the case of  periodic domains by L. Carbone and F. Colombini \cite{CC} and then in a more general framework by E. De Giorgi, G. Dal Maso and P. Longo \cite{DDL} and G. Dal Maso and P. Longo \cite{DL}, G. Dal Maso \cite{D}.
Our main reference for this work will be the papers of D. Cioranescu and F. Murat \cite{CM1,CM2}, in which the case of  a periodic repartition of the holes $D\setminus\De$ is studied. It is  proved that when the number of holes and their size  are evolving in a critical fashion, then the limiting problem is no longer an obstacle problem, but a simple elliptic boundary value problem with a new term that takes into account the effect of the holes. 

Our goal is to generalize their result to the case where 
the holes are still located in small neighborhoods of the points of the lattice $\eps\ZZ^{n}$ but have random size and shape. 
More precisely, we assume that for any $\eps$ and $\omega$ the domain
$\De(\omega)$  is obtained from a fixed set $D$ by perforating holes $\{S_{\eps}(k,\omega)\, ;\,  k\in \ZZ^{n}\}$ such that
$$
S_{\eps} (k,\omega) \subset  B_{\eps/2}(\eps k)\qquad  \mbox{ for all } k\in \ZZ^{n}.
$$  
We denote by 
$$ T_{\eps}(\omega)= \cup_{k\in\ZZ^{n}} S_{\eps}(k,\omega) \cap D$$
the union of all the holes in $D$. We then have
$$ \De (\omega )= D(\omega)\setminus T_{\eps}(\omega).$$
The assumptions on the sets $S_{\eps} (k,\omega)$ will be made precise in the next section. We can already point out the fact that we will not exclude the case where $S_{\eps}(k,\omega)=\emptyset$ for some $k$, thus allowing the fact that no holes may be present at some lattice points.
\vspace{10pt}

With these notations, we rewrite the obstacle problem as follows:
\begin{equation} \label{eq:inf}
\J (u^\eps) = \inf_{v \in K^\eps } \J (v) ,\quad \qquad u^\eps\in K_\eps
\end{equation}
with
$$\J(v) = \int_D \frac{1}{2} |\na u|^2 - f \, u\, dx$$
and
$$K_\eps = \{v\in H^1_0 (D)\, ;\, v \geq 0 \mbox{ a.e. in } T_\eps\}.$$
Since $K_{\eps}$ is closed, convex and not empty, (\ref{eq:inf}) has a unique solution $u^{\eps}\in K_\eps$.
Moreover, $u^\eps$ solves
\begin{equation} \label{eq:0}
\left\{
\begin{array}{ll}
 - \Delta u^{\eps} = f & \mbox{ in } \De \\[5pt]
u^{\eps}(x) \geq 0  & \mbox{ on } T_{\eps} \\[5pt]
u^{\eps} (x)= 0 & \mbox{ on }  \pa  D \setminus T_{\eps}
\end{array}
\right.
\end{equation}
As mentioned in the introduction, it is expected that under appropriate assumptions on the size of the holes $S_\eps(k,\omega)$, the function $u^\eps$ converges weakly in $H^1$ to $u$ solution of
$$
\begin{array}{ll}
- \Delta u - \alpha_0 u_- = f & \mbox{ in } D \\[5pt]
u=0 & \mbox{ on } \pa D.
\end{array}
$$
where $u_-(x)=max(0,-u(x))$.
\vspace{10pt}

The assumptions and the result are made precise in the next section.
The proof of the main theorem, which is details in Section \ref{sec:thm}, relies on the construction of an appropriate corrector. This construction is detailed in Sections \ref{sec:balls} and \ref{sec:general}, first in the case where the holes are balls in dimension $n\geq 2$, then when no assumptions are made on the shape of the holes (in dimension $n\geq 3$ only).

\section{Assumptions and Main result}
First, we need  to make precise our assumptions on the holes $S_\eps(k,\omega)$. The first assumption is mainly technical:
\vspace{10pt}

{\bf Assumption 1:}  There exists a (large) constant $M$ such that for all $k\in \ZZ^n$ and a.e. $\omega\in \Omega$ we have
 $$
 \begin{array}{ll}
 S_\eps(k,\omega)\subset B_{M \eps^{n/(n-2)}}(\eps k) & \mbox{ if } n\geq 3\\[5pt]
 S_\eps(k,\omega)\subset B_{\exp(-M \eps^{-2})}(\eps k) & \mbox{ if } n=2 
\end{array}
$$
for $\eps$ small.

As mentioned in the introduction, the asymptotic behavior of the $u^\eps$ strongly depends on the size of the holes. 
The critical size for which interesting phenomena is observed corresponds to finite, non trivial capacity of the set $T_\eps$. More precisely, we assume: 
\vspace{10pt}

{\bf Assumption 2:} 
For all $k\in \ZZ^n$ and a.e. $\omega\in \Omega$, there exists $\gamma(k,\omega)$ (independent of $\eps$) such that
$$ \capa (S_\eps(k,\omega)) = \eps^n \gamma(k,\omega),$$
where $\capa (A)$ denote the capacity of subset $A$ of $\RR^n$, defined by:
$$ \capa (A) = \inf\left\{\int_{\RR^n} |\na h|^2\, dx\, ;\; h\in H^1(\RR^n),\; h\geq 1 \mbox{ in } A, \lim_{|x|\rightarrow \infty} h(x)= 0\right\},$$
in dimension $n\geq 3$ and by 
$$ \capa (A) = \inf\left\{\int_{B_1} |\na h|^2\, dx\, ;\; h\in H^1_0(B_1),\; h\geq 1 \mbox{ in } A \right\},$$
in dimension $n=2$ and for sets $A\subset B_1$.
Moreover, we assume that there exists a constant $\overline \gamma>0$:
\begin{equation}\label{eq:gb} 
 \gamma (k,\omega) \leq \overline \gamma \quad \mbox{ for all $k\in\ZZ^n$ and a.e. $\omega\in\Omega$}.
 \end{equation}

\vspace{10pt}

Finally, our last assumption will be necessary to ensure that some averaging process occur as $\eps$ goes to zero:
\vspace{10pt}

{\bf Assumption 3:} 
The process
$ \gamma : \ZZ^{n} \times \Omega \mapsto [0,\infty)$
is stationary ergodic:  There exists a family of measure-preserving transformations $\tau_k: \Omega\rightarrow \Omega$ satisfying
$$ \gamma(k+k',\omega)=\gamma(k,\tau_{k'}\omega)\quad \mbox{ for all } k,k'\in\ZZ^n\mbox{ and } \omega\in\Omega,$$
and such that 
if $A\subset \Omega$ and  $\tau_k A=A $ for all $k\in\ZZ^n$, then $P(A)=0$ or $P(A)=1$
(the only invariant set of positive measure is the whole set).

\vspace{15pt}

Let us make a few remarks concerning those assumptions:
First of all, we stress out the fact that the shape of the holes $S_\eps$ is left unspecified and may change with $\eps$; Only the rescaled capacity is independent on $\eps$. 
The first assumption, which implies that the diameters of the holes decrease faster than $\eps$,  guarantees that the capacities of neighboring sets separate at the limit (i.e. that $\capa(\cup S_\eps) \sim \sum \capa(S_\eps)$). And the choice of scaling for the capacity guarantee that $\capa(T_\eps) $ remains bounded as $\eps$ goes to zero (since $\#\{ \ZZ^n\cap \eps^{-1} D\} \leq C\eps^{-n}$).
Finally, the hypothesis of stationarity is the most general extension of the notions of periodicity and almost periodicity for a function to have some self-averaging behavior.

\vspace{15pt}

Under those assumptions, we prove the following result:
\begin{theorem}\label{thm:1}
Assume that $n\geq 3$ or $n=2$ and the holes are all balls. Then
there exists $\alpha_0\geq 0$ such that when  $\eps$ goes to zero,
$u^\eps$ converges weakly in $H^1$ to a function $\ou$  solution of
the following minimization problem
$$ \min \left\{  \int_{D}  \frac{1}{2} |\na u|^2 +\frac{1}{2} \alpha_0 u_-^2 - f \, u\, dx\, ;\,  u \in  H^1_0 (D) \right\},
$$
where $u_-(x)=\max(0,-u(x)).$
In particular, $\ou$ solves
$$\left\{
\begin{array}{ll}
- \Delta \ou - \alpha_0 \ou_- = f & \mbox{ in } D \\[5pt]
\ou=0 & \mbox{ on } \pa D.
\end{array}
\right.
$$
Moreover, if there exists $\underline \gamma>0$ such that 
$$\gamma (k,\omega) \geq \underline \gamma \quad \mbox{ for all $k\in\ZZ^n$ and a.e. $\omega\in\Omega$},$$
then $\alpha_0>0$.
\end{theorem}

The general result holds also in dimension $n=2$ when the holes have random shape. However, because  the fundamental solution of Laplace's equation is different in that case, the proof is slightly different and more technical.
\vspace{10pt}

As in Cioranescu - Murat \cite{CM1,CM2}, the proof of this result relies on the construction of an appropriate corrector.
More precisely, the key is the following result:
\begin{proposition}\label{prop:1}
Under the assumptions listed above, 
there exists  a non-negative real number $\alpha_0$ and a function $w^\eps_0(x,\omega)$ such that 
$$
\left\{
\begin{array}{l}
 \Delta w^\eps  = \alpha_0 \quad \mbox{ in }  \De(\omega) \\[5pt]
 w^\eps (x) = 1\quad \mbox{ in } T_\eps(\omega)\\[5pt]
 w^\eps (x) = 0\quad \mbox{ on } \pa D \setminus T_\eps(\omega) 
 \end{array}
\right.
$$
for almost all $\omega \in \Omega$, and 
$$  w^\eps  \longrightarrow 0\qquad H^1(D)\mbox{-weak}\mbox{ a.s. } \omega \in \Omega.$$
\end{proposition}

Note that as in \cite{CM1}, the equation 
$$ \Delta w^\eps  = \alpha_0 \quad \mbox{ in }  \De(\omega) $$
can be replaced by the weaker condition:
\begin{equation}\label{H1'}
\left\{
\begin{array}{l}
\mbox{For all sequences $v^\eps$ satisfying:} \\[5pt]
\quad\qquad\left\{ \begin{array}{l}
v^\eps = 0\quad \mbox{ on } T_\eps\\
v^\eps\longrightarrow v\quad \mbox{ in } H^1(D)-\mbox{weak}
\end{array}\right.\\[10pt]
\mbox{and for any $\phi\in\mathcal D(D)$, we have:}\\[5pt]
\qquad \langle \Delta w^\eps \, ,\,\phi v^\eps \rangle_{H^{-1},H^1_0(D)} \longrightarrow \langle \alpha_0\, , \, \phi v\rangle.
\end{array}
\right.
\end{equation}

The proof of Proposition \ref{prop:1}  will occupy most of this paper.
It will be split in two parts: In Section \ref{sec:balls}, we consider the (simpler) case when the holes $S_\eps(k,\omega)$ are all balls of random radius. In Section \ref{sec:general}, we will use this first result to treat the general case (when the holes have unspecified shapes).

Before turning to this proof, we briefly give, in the next section the proof of the main theorem.
\vspace{25pt}

\section{Proof of Theorem \ref{thm:1}} \label{sec:thm}
First of all, standard elliptic estimates give the existence of a function $\ou$ such that  
$$ u^\eps \longrightarrow \ou \qquad H^1-\mbox{weak}.$$ 
If we introduce the limit energy
$$ \Ja (v) =\int_D \frac{1}{2} |\na v|^2 +\frac{1}{2} \alpha_0 v_-^2- f \, v\, dx,$$
it is readily seen that all we need to show is the following inequality:
$$ \Ja (\ou) = \inf_{v\in H^1_0(D)} \Ja(v),$$
\vspace{10pt}

This relies on the following two lemmas:
\begin{lemma}\label{lem:1}
For any $\vphi \in W^{2,\infty}_0$, we have
$$ \lim_{\eps\rightarrow 0} \int_D  |\na w^\eps |^2 \vphi\, dx = \int _D \alpha_0 \vphi\, dx$$
\end{lemma}

\begin{lemma}\label{lem:2}
If $u^\eps \rightharpoonup \ou$ in  $H^1$-weak, then
$$ \liminf_{\eps\rightarrow 0} \J(u^\eps) \geq \Ja(\ou)$$
\end{lemma}

Let us see that those two lemmas imply the theorem:
For any $v\in W^{1,\infty}_0$, the function $v+v_- w^\eps $ is non-negative on the holes, and is thus admissible for the initial obstacle problem.
In particular by definition of $u^\eps$, we have
$$ \J(u^\eps) \leq \J(v+v_- w^\eps ).$$
We write
\begin{eqnarray*}
\J(v+v_- w^\eps ) & = &  \int \frac{1}{2} \big[ |\na v|^2 + |\na v_-|^2 {w^\eps }^2+|v_-|^2|\na w^\eps |^2\big]\, dx \\
&& +\int \big[ v_-\na v_- w^\eps \na w^\eps + \na v\na v_- {w^\eps }+\na v v_-\na w^\eps \big]\, dx
\end{eqnarray*}
and it is readily check that Lemma \ref{lem:1} and the weak convergence of $w^\eps $ to $0$ in $H^1$ implies 
$$ \lim _{\eps\rightarrow 0}  \J (v+v_- w^\eps ) = \Ja (v),$$
as soon as $v\in W^{2,\infty}$.
We deduce:
$$ \Ja(v) \geq \limsup_{\eps\rightarrow 0} \J(u^\eps)$$
for all $v\in W^{2,\infty}_0$.
Together with Lemma \ref{lem:2} this gives
$$ \Ja(v) \geq \Ja(\ou)$$
for all $v\in W^{2,\infty}_0$.
We deduce Theorem \ref{thm:1} by a density argument.
\qed

\vspace{10pt}

{\bf Proof of Lemma \ref{lem:1}:} We recall the proof of Cioranescu-Murat \cite{CM1}: Since $1-w^\eps =0$ in $T_\eps$, we have:
\begin{eqnarray*} 
\int_{\De} \Delta w^\eps  \vphi (1-w^\eps )\, dx 
&=& \int_{D_\eps} \vphi |\na w^\eps |^2\, dx - \int_{D_\eps}\na \vphi \na w^\eps   (1-w^\eps )\, dx
\end{eqnarray*}
and so
\begin{eqnarray*} 
\int_{D}\alpha_0 \vphi (1-w^\eps )\, dx
& =&  \int_{D_\eps}\alpha_0 \vphi (1-w^\eps )\, dx \\
& =&  \int_{D_\eps} \vphi |\na w^\eps |^2\, dx - \int_{D_\eps}\na \vphi \na w^\eps   (1-w^\eps )\, dx.
\end{eqnarray*}
finally
\begin{eqnarray*}
\int_{D_\eps}\na \vphi \na w^\eps   (1-w^\eps )\, dx & = &  \int_{D_\eps}\na \vphi \na w^\eps   -\int_{D_\eps} \na \vphi \na w^\eps  w^\eps \, dx \longrightarrow  0
\end{eqnarray*}
since $w^\eps $ goes to zero $H^1$-weak and $L^2$-strong.
The lemma follows.
\qed

 \vspace{10pt}

{\bf Proof of Lemma \ref{lem:2}:} See  Cioranescu-Murat \cite{CM2}, Proposition~3.1.

\vspace{20pt}

\section{Proof of Proposition \ref{prop:1}: Balls of random radius} \label{sec:balls}

Throughout this section, we assume that the sets $S_\eps(k,\omega)$ are balls centered at $\eps k$.
Since
$$ \capa (B_r)=\left\{ \begin{array}{ll}
n(n-2)\omega_n r^{n-2} &\mbox{ if } n\geq 3, \\[8pt]
\ds -\frac{2\pi}{\log r} & \mbox{ if } n=2
\end{array}
\right.
$$ 
Assumption 2 becomes in this framework:
$$
 S_{\eps}(k,\omega) = 
 B_{a^\eps(r(k,\omega)) }(\eps k) \quad \mbox{ for all } k\in\ZZ^{n}
$$
with
$$
a^\eps(r)=
\left\{
\begin{array}{ll}
r \eps^{n/(n-2)}  & \mbox{ if } n\geq 3, \\[5pt]
 \exp(-r^{-1}\eps^{-2}) & \mbox{ if } n=2 ,
  \end{array}
  \right. 
 $$
and
$$
r(k,\omega) 
=
\left\{
\begin{array}{ll}
\left(\frac{\gamma(k,\omega)}{n(n-2)\omega_n}\right)^{1/(n-2)} & \mbox{ if } n\geq 3, \\[5pt]
 \gamma(k,\omega)/2\pi  & \mbox{ if } n=2 .
  \end{array}
  \right. 
 $$
Note in particular that the process
$$ r : \ZZ^{n} \times \Omega \mapsto [0,\infty)$$
is stationary ergodic and satisfies
\begin{equation}\label{eq:rb} 
 r(k,\omega) \leq \overline r \quad \mbox{ for all $k\in\ZZ^n$ and a.e. $\omega\in\Omega$}
 \end{equation}
for some constant $\overline r>0$.
Without loss of generality, we can always assume that $\overline r < 1/2$ (so that there is no overlapping of the holes for any $\eps<1$):

\vspace{10pt}

\subsection{The auxiliary obstacle problem }
After rescaling, we  look for the corrector $w^\eps(x,\omega)$ in the form
$$w^\eps(x,\omega)=\eps^{2} v^\eps(x/\eps,\omega)$$
with $v^\eps(y,\omega)$ solution to 
$$ \left\{
\begin{array}{ll}
\ds \Delta v = \alpha, \qquad& \mbox{ in}  \eps^{-1} \De\,, \quad \mbox{a.e. } \omega\in \Omega\\ 
\ds v = \eps^{-2}& \mbox{ on } \cup_{k\in\ZZ^{n}}B_{ \overline a^\eps (k,\omega) } (k)
\end{array}
\right.
$$
with 
$$
\overline a^\eps(r)= \left\{
\begin{array}{ll}
r \eps^{2/(n-2)} & \mbox{ if } n\geq 3, \\[5pt]
 \eps^{-1} \exp(-r^{-1}\eps^{-2}) & \mbox{ if } n=2 ,
  \end{array}
  \right.  $$
and satisfying
$$  \eps^2 v^\eps(x/\eps)\longrightarrow 0 \mbox{ in $H^1$-weak }.$$
\vspace{10pt}

One of the main tool in the proof is the fundamental solution of the Laplace equation, given by:
$$
h(x)=\left\{
\begin{array}{ll}
\ds \frac{1}{n(n-2)\omega_{n}}\frac{1}{|x|^{n-2}} & \mbox{ if } n\geq 3, \\[8pt]
\ds -\frac{1}{2\pi} \log|x| & \mbox{ if } n=2.
\end{array}
\right.
$$
In particular, we note that
$$ h|_{B_{\overline a^\eps(r(k,\omega)) } (0)} = 
\left\{
\begin{array}{ll}
\ds  \frac{1}{n(n-2)\omega_{n}r^{n-2} } \eps^{-2} & \mbox{ if } n\geq 3, \\[8pt]
\ds \frac{1}{2\pi}(\log(\eps)+r^{-1}\eps^{-2}) & \mbox{ if } n=2,
 \end{array}\right.
$$
so we expect the rescaled  corrector $v^\eps(x,\omega)$ to behave near the hole $B_{\overline a (k,\omega)}(k)$  like the function 
$$h_k(x):=
\left\{
\begin{array}{ll}
\ds
 \gamma(k,\omega)\;  h(x-k) = \frac{r(k,\omega)^{n-2}}{|x-k|^{n-2}},& \mbox{ if } n\geq 3 \\[8pt]
\ds  \gamma(k,\omega)\; h(x-k)   = -r(k,\omega) \log |x-k|& \mbox{ if } n=2,
 \end{array}\right.
$$
where
$$\gamma(k,\omega) =\left\{
\begin{array}{ll}
\ds  (r(k,\omega))^{n-2}n(n-2)\omega_{n}& \mbox{ if } n\geq 3 \\[5pt]
\ds 2\pi r(k,\omega) & \mbox{ if } n=2.
 \end{array}\right.
$$

Since $h_k$ satisfies
$$  \Delta h_k =  - \gamma(k,\omega)\; \delta(x-k),$$
we will construct $v^\eps(x,\omega)$ by solving
$$
\left\{\begin{array}{l}
 \Delta v = \alpha - \sum_{k\in\ZZ^{n}\cap A} \gamma(k,\omega) \delta(x-k) \,\quad \mbox{ in } D\, , \qquad \\[5pt]
 v=0  \quad\mbox{ on } \pa D.
\end{array}
\right.
$$
The main issue is thus to find the critical $\alpha$ for which the solution of the above equation has the appropriate behavior near $x=k$.

Following \cite{CSW}, this will be done by  introducing the following obstacle problem, for every open set $A\subset \RR^n$ and $\alpha\in\RR$:
\begin{equation}\label{def:ov}
\ov_{\alpha,A}(x,\omega) = \inf\left\{v(x)\, ; \, \Delta v \leq \alpha - \sum_{k\in\ZZ^{n}\cap A} \gamma(k,\omega) \delta(x-k) \, ,\,
\begin{array}{l} 
v\geq 0  \mbox{ in } A \\
 v=0  \mbox{ on } \pa A
 \end{array}
 \right\}.
\end{equation}

Clearly, the function $\ov_{\alpha,A}$ is solution of 
\begin{equation} \label{eq:v}
\Delta v = \alpha - \sum_{k\in\ZZ^{n}\cap A} \gamma(k,\omega) \delta(x-k) 
\end{equation}
whenever it is positive.
Note that the function
\begin{eqnarray} \label{eq:hak} 
h_{\alpha,k}(x) & : = & \frac{\alpha}{2n} |x-k|^2 + h_k(x-k)  \\[5pt]
& = & \left\{
\begin{array}{ll}
\ds
\frac{\alpha}{2n} |x-k|^2+ \frac{r(k,\omega)^{n-2}}{|x-k|^{n-2}},& \mbox{ if } n\geq 3, \\[8pt]
\ds\frac{\alpha}{2n} |x-k|^2 -r(k,\omega) \log |x-k|& \mbox{ if } n=2,
 \end{array}\right. \nonumber
\end{eqnarray}
also satisfies
$$ \Delta h_{\alpha,k} (x) = \alpha - \gamma(k,\omega) \delta(x-k).$$
It follows from  (\ref{eq:v})  and the maximum principle that 
if $B_1(k)\subset A$, then, for all $x$ in $B_1(k)$ and for  almost every $\omega $ in $\Omega$, we have
\begin{equation}\label{eq:ineqw} 
\ov_{\alpha,A}(x,\omega)  \geq 
\left\{
\begin{array}{ll}
\ds h_{\alpha,k}(x) - \frac{\alpha}{2n} - r^{n-2} & \mbox{ if } n\geq 3 \\[8pt]
\ds h_{\alpha,k}(x)- \frac{\alpha}{2n} & \mbox{ if } n=2.
\end{array}\right.
\end{equation}

\vspace{10pt}

\subsection{Critical $\alpha$}
The purpose of this section is to prove that for a critical $\alpha$, $\ov_{\alpha,A}$ behaves like $h_{\alpha,k}$ near $S_\eps(k,\omega)$.
For that purpose, we introduce the following quantity, which measures the size of the contact set: 
$$
\om_{\alpha}(A,\omega) = |\{x\in A\, ;\, \ov _{\alpha,A}(x,\omega)= 0 \}| 
$$
where $|A|$ denotes the Lebesgue measure of a set $A$.
\vspace{10pt}

The starting point of the proof is the following lemma:
\begin{lemma}\label{lem:ergodic}
The random variable $\om_{\alpha}$ is subadditive, 
and the process 
$$T_k m(A,\omega) = m(k+A,\omega)$$
has the same distribution for all $k\in \ZZ^n$.
\end{lemma}
{\bf Proof of  Lemma \ref{lem:ergodic}:}
Assume that the finite family  of sets $(A_i)_{i\in I}$ is such that 
$$ 
\begin{array}{l}
A_i \subset A \qquad \mbox{ for all } i\in I \\
A_i\cap A_j = \emptyset \quad  \mbox{ for all } i\neq j \\
|A-\cup_{i\in I} A_i| = 0
\end{array}
$$ 
then $\ov_{\alpha,A}$ is admissible for each $A_i$, and so $\ov_{\alpha,A_i} \leq u_{\alpha,A}$.
It follows that 
$$ \{\ov_{\alpha,A}=0\}\cap A_i \subset\{\ov_{\alpha,A_i}=0\}$$
and so
$$
\om_{\alpha}(A,\omega)  = \sum_{i\in I} |\{\ov_{\alpha,A}=0\}\cap A_i | \leq  \sum_{i\in I}|\{\ov_{\alpha,A_i}=0\}| =   \sum_{i\in I}\om_{\alpha}(A_i,\omega),$$
which gives the subadditive property.
Assumption {\bf 3} then yields
$$ T_k m(A,\omega) = m(A,\tau_k \omega)$$
which gives the last assertion of the lemma.
\qed
\vspace{20pt}

Since 
$ \om_{\alpha}(A,\omega) \leq |A|,$ and thanks to the ergodicity of the transformations $\tau_k$,
it follows from the subadditive ergodic theorem (see \cite{DM}) that for each $\alpha$, there exists a constant $\oell(\alpha)$ such that
$$ \lim_{t\rightarrow \infty} \frac{\om_{\alpha}( B_t(0),\omega)}{|B_t(0)|} = \oell (\alpha)\quad  \mbox{ a.s., } \quad $$
where $B_t(0)$ denotes the ball centered at the origin with radius $t$.
Note that the limit exists and is the same if instead of $B_t(0)$, we use cubes or balls  centered at $tx_0$ for some $x_0$.

If we scale back and consider the function 
$$ \overline  w^\eps_\alpha (y,\omega)= \eps^2 \; \ov_{\alpha,B_{\eps^{-1}}(\eps^{-1}x_0)}(y/\eps,\omega) ,\qquad \mbox{ in } B_1(x_0),$$
we deduce
$$ \lim_{\eps\rightarrow 0} \frac{|\{ y\, ;\, \overline w^\eps_{\alpha}(y,\omega)=0\}|}{|B_1|} = \oell (\alpha)\quad  \mbox{ a.s. }
$$

\vspace{10pt}

The next lemma summarizes the properties of  $\oell(\alpha)$:
\begin{lemma}\label{lem:alpha}
\item[(i)] $\oell(\alpha)$ is a nondecreasing functions of $\alpha$.
\item[(ii)] If $\alpha<0$, then $\oell (\alpha) = 0$. Moreover, if the radii $r(k,\omega)$ are bounded from below, then $\oell (\alpha) = 0$ for any $\alpha $  such that $\alpha <n(n-2) \inf_{k\in\ZZ^n} r(k,\omega)^{n-2}  $ almost surely.
\item[(iii)] If $\alpha\geq  2^n n(n-2)\sup_{k\in\ZZ^n} r(k,\omega)^{n-2}$ (or $\alpha \geq 8 r$ for $n=2$) almost surely, then
 $\oell (\alpha) > 0$.
\end{lemma}
{\bf  Proof.} \\
\noindent (i) The proof follows immediately from the inequality
$$ \ov_{\alpha,A} \leq \ov_{\alpha',A}\qquad  \mbox{ for any $\alpha$, $\alpha'$ such that } \alpha' \leq \alpha.$$

\noindent (ii)
If $\alpha$ is negative, then the function $\frac{\alpha}{2n}|x-x_0|^2-\frac{\alpha}{2n} (tr)^2$, which is a  sub-solution of (\ref{eq:v}), is positive in $t B_r(x_0)$ and  vanishes along $\pa (t B_r(x_0))$ for any ball $B_r(x_0)$ and for any $t>0$. We deduce:
$$ \ov_{\alpha,tB} > \frac{\alpha}{2n}|x-x_0|^2-\frac{\alpha}{2n} (tr)^2 >0 \mbox{ in } t B_r(x_0)$$
for all $t>0$.  Therefore $m_\alpha (tB,\omega) =0$ for all $t>0$, so $\oell(\alpha)=0$ for all $\alpha<0$.

Furthermore, if $r(k,\omega)$ is bounded below:
$$\;\; r(k,\omega) \geq \underline r >0 \mbox{ for all $k\in\ZZ^n$, a.e. $\omega\in\Omega$}, $$
then, the function
$\frac{\alpha}{2n} |x-k|^2+  \frac{ \underline r^{n-2}}{|x-k|^{n-2}}-\frac{\alpha}{2n}-\underline r^{n-2}$
is a solution of (\ref{eq:v}) in $B_1(k)$ which vanishes on $\pa B_1(k)$ and  is strictly positive in $B_1(k)$ as long as $\alpha<n(n-2)\underline r^{n-2}$.
As above, we deduce that   $m_\alpha (tB,\omega) =0$ for all $t>0$ and for all $ \alpha<n(n-2)\underline r^{n-2}$.
\vspace{10pt}

\noindent (iii) 
The function $h_{\alpha,k}(x)=\frac{\alpha}{2n}|x-k|^2+\frac{r^{n-2}}{|x-k|^{n-2}}$  is radially symmetric and reaches its minimum when 
\begin{equation}\label{eq:Rak}|x-k| =R(\alpha,k):= 
\left\{
\begin{array}{ll}
\left(\frac{n(n-2)r(k,\omega)^{n-2}}{\alpha}\right)^{1/n} &\mbox{ when } n\geq 3 \\[5pt]
\left(\frac{2 r(k,\omega)}{\alpha} \right)^{1/2} & \mbox{ when } n=2 
\end{array}
\right.
\end{equation}

In particular, for  $\alpha> 2^n n(n-2) r(k,\omega)^{n-2}$ (or $n\geq 8r(k,\omega$ when $n=2$), we have $R(\alpha,k) <1/2$ and so the function
$$g_k(x)=
\left\{
\begin{array}{ll}
h_{\alpha ,k}(x)-  D_k & \mbox{ in } B_{R(\alpha,k)}(k)\\
0 &\mbox{  in }  \RR^n\setminus B_{R(\alpha,k)}(k)
\end{array}
\right.
$$
satisfies
$$
\Delta g_k \leq \alpha - \gamma(k,\omega)\delta(x-k) \mbox{ in } C_1(k)\, ,
$$
and
$$ 
g_k=0\quad \mbox{ in } C_1(k) \setminus B_{1/2}(k)
$$
where $C_1(k)$ denotes the cube of size $1$ centered at $k$, and the constant $C_k$ is chosen in such a way that $g_k$ and $\na g_k$ vanish along $\pa B_{R(\alpha,k)}$:
\begin{equation}\label{eq:Dak}
D(\alpha,k):=\left\{
\begin{array}{ll}
 (\frac{\alpha}{2n})^{\frac{n-2}{n}}r^{\frac{2(n-2)}{n}} \left(\frac{n-2}{2}\right)^{\frac{2}{n}}\left(\frac{n}{n-2}\right)  &\mbox{ when } n\geq 3 \\[5pt]
\frac{r}{2}\left(1-\log(2r/\alpha) \right) & \mbox{ when } n=2 .
\end{array}
\right.
  \end{equation}

By definition of $\ov_{\alpha,tB}$, we deduce that
$$ \ov_{\alpha,tB}(x) \leq \sum_{k\in\ZZ^n\cap tB} g_k (x)\quad \mbox{ in } t B \mbox{ a.s.} $$
In particular, this implies that $ \ov_{\alpha,tB}$ vanishes in $t B \setminus \cup_{k\in\ZZ^n} B_{1/2}(k)$,
and so 
$$ \frac{\overline m_\alpha(t B,\omega) }{|tB|}\geq  \left(\frac{|C_1| -|B_{1/2}|}{|C_1|}\right)  = 1-\frac{\omega_n}{2^n}\quad\mbox{ a.s. }$$  
We conclude
$$ \oell(\alpha ) \geq 1-\frac{\omega_n}{2^n}>0.
$$ \qed

\vspace{20pt}

Using Lemma \ref{lem:alpha},  we can define
$$\alpha_{0} = \sup \{\alpha\, ;\, \oell(\alpha) = 0 \}.$$
Note that $\alpha_0$ is finite under Assumption 3 (Lemma \ref{lem:alpha} (iii)) and that $\alpha_0\geq 0$ is strictly positive as soon as the $r(k,\omega)$ are bounded from below almost surely by a positive constant (Lemma \ref{lem:alpha} (ii)).

In the rest of this section, we are going to show that the function
$$
 w^\eps (x,\omega) = \inf\left\{w(x)\, ; \, \Delta w \leq \alpha_0 \mbox{ in } D\setminus T_\eps \, ,\,
\begin{array}{l} 
w \geq 1  \mbox{ on } T_\eps \cap D \\
w = 0  \mbox{ on } \pa D \setminus T_\eps 
 \end{array}
 \right\},
$$
satisfies all the conditions of Proposition \ref{prop:1}. 
We will rely on a series of intermediate functions. 

For the first lemma, we fix a bounded subset $A$ of $\RR^n$ and 
we  denote by 
\begin{equation}\label{eq:ovae}
\ovae(x,\omega) =\ov_{\alpha,\eps^{-1}A}(x,\omega)
\end{equation}
the solutions of (\ref{def:ov}) defined in $\eps^{-1}A$. 
We also introduce  the rescaled function
$$ \wae (y,\omega)= \eps^2 \; \ovae (y/\eps,\omega) ,$$
defined in $A$.
\vspace{20pt}

The key properties of $\ovae $ are given by the following lemma:
\begin{lemma}\label{lem:1.2}
\item[(i)] For every $\alpha $ and for every $k\in\ZZ^n$, we have 
$$ \ovae(x)\geq
\left\{
\begin{array}{ll}
\ds h_{\alpha,k}(x) - \frac{\alpha}{2n} - r^{n-2} & \mbox{ if } n\geq 3 \\[8pt]
\ds h_{\alpha,k}(x)- \frac{\alpha}{2n} & \mbox{ if } n=2
\end{array}\right.
$$
for all $x\in B_{1}(k)$ and  almost everywhere $\omega\in\Omega$ (where $ h_{\alpha,k}$ is defined by (\ref{eq:hak})).
\item[(ii)] For every $\alpha > \alpha_0$, we have 
$$ \ovae(x)\leq h_{\alpha,k}(x)+ o(\eps^{-2})$$ 
for all $x\in B_{1/2}(k)$ and  almost everywhere $\omega\in\Omega$.
\end{lemma}
Since
$$ h_{\alpha,k}|_{B_{\overline a^\eps(r(k,\omega)) } (0)} = 
\left\{
\begin{array}{ll}
\ds  \eps^{-2} +\frac{ \alpha_0 }{2n}|\overline a^\eps(r(k,\omega)) |^2 & \mbox{ if } n\geq 3 \\[8pt]
\ds \eps^{-2}+  \frac{ \alpha_0 }{4}|\overline a^\eps(r(k,\omega)) |^2 +r(k,\omega)\log \eps  & \mbox{ if } n=2,
 \end{array}\right.
$$
we deduce the following corollary:
\begin{corollary}
\item[(i)] For every $\alpha $ and every $k\in\ZZ^n$ such that $r(k,\omega)>0$, we have 
$$ \ovae(x)\geq \eps^{-2}+o(1) \quad \mbox{ on } \pa B_{\overline a^\eps(r(k,\omega))}(k) \quad\mbox{ a.e. } \omega\in\Omega$$
and so
$$ \wae \geq 1+ o(\eps^2) \mbox{ on } \pa T_\eps(\omega) \quad\mbox{ a.e. } \omega\in\Omega$$
for all $\alpha$. 
\item[(ii)] For every $\alpha > \alpha_0$ and every $k\in\ZZ^n$, we have 
$$ \ovae(x)\leq \eps^{-2}+o(\eps^{-2})\quad \mbox{ on } \pa B_{\overline a^\eps(r(k,\omega))}(k) \quad\mbox{ a.e. } \omega\in\Omega$$
and so
$$ \wae \leq 1+ o(1) \mbox{ on } \pa T_\eps(\omega)\quad\mbox{ a.e. } \omega\in\Omega$$

\end{corollary}
{\it Proof of Lemma \ref{lem:1.2}:}\\
\noindent (i) Immediate consequence of (\ref{eq:ineqw}).
\vspace{10pt}

\noindent (ii) \\
{\it Preliminary:}
First of all since $A$ is bounded, we have
$$ A \subset B_R(x_0).$$
Without loss of generality, we can  always assume that $B_R(x_0)=B_1(0)$.
We then introduce
$$
\vae(x,\omega) =\ov_{\alpha,\eps^{-1}B_1}(x,\omega),
$$
the solutions of (\ref{def:ov}) in $B_{\eps^{-1}}(0)$.
It is readily seen that $\vae$ is admissible for (\ref{def:ov}) and thus
$$\ovae(x,\omega) \leq \vae(x,\omega)  \qquad \mbox{for all } x\in \eps^{-1}A\, \; \mbox{ a.e. }\omega\in\Omega.$$
It is thus enough to prove (ii) for $\vae$.
\vspace{10pt}

We will need the following consequence of Lemma \ref{lem:ergodic} (see \cite{CSW} for the proof):
\begin{lemma}\label{lem:local}
For any  ball $B_r(x_0)\in B_1(0)$, the following limit holds, a.s. in $\omega$
$$\lim_{\eps\rightarrow 0} \frac{| \{\vae (x,\omega)= 0 \}\cap B_{\eps^{-1}r}(\eps^{-1}x_1)  | }{|B_{\eps^{-1}r}|} = \oell(\alpha)
$$
\end{lemma}
\vspace{10pt}

\noindent {\it Step 1:} We can now start the proof:
For any $\delta>0$, we can cover $B_{\eps^{-1}}$ by a finite number $N$ ($\leq C\delta^{-n}$) of balls $B_i$ with radius $\delta \eps^{-1}$ and center $\eps^{-1}x_i$.
Since $\alpha >\alpha_0$, we have $\oell(\alpha)>0$. By Lemma \ref{lem:local}, we deduce that for every $i$,  there exists $\eps_i$ such that if $\eps\leq\eps_i$, then
$$| \{\vae (x,\omega)= 0 \}\cap B_{i}  | >0\quad \mbox{ a.s. } \omega.$$
In particular, if $\eps \leq \inf \eps_i$, then $\vae(y_i)=0$ for some $y_i$ in $B_i$ a.s. $\omega\in\Omega$.
We now have to show that this implies that  $\vae$ remains small in each $B_i$ as long as we stay away from the lattice points $k\in\ZZ^n$.
More precisely, we want to show that
$$ \sup_{B_i\setminus\cup_{k\in\ZZ^n} B_{1/4}(k)} \vae \leq C\delta^2 \eps^{-2}.$$

\noindent {\it Step 2:} Let $\eta$ be a nonnegative function such that $0\leq \eta(x)\leq 1$ for all $x$,  $\eta(x)=1$ in $B_{1/8}$ and $\eta = 0$ in $\RR^n\setminus B_{1/4}$.
Then the function $u=\vae \star \eta$ is nonnegative on $2B_i$ and satisfies 
$$ -C\leq \Delta u \leq C$$
where $C$ is a universal constant depending only on $n$ and $\overline r$.
In particular, since $B_i$ has radius $\delta \eps^{-1}$, Harnack inequality yields: 
$$\sup_{B_i} u \leq C \inf_{B_i}   u + C\alpha (\delta\eps^{-1})^2.$$
\vspace{10pt}

\noindent {\it Step 3:} We need the following lemma:
\begin{lemma}\label{lem:average}
If $\Delta v \leq \alpha$ in $B_r(y_0)$, then
$$ \frac{1}{B_r}\int_{B_r(y_0)} v (x)\, dx \leq v (y_0) + \alpha C(n) r^2$$
where $C(n)$ is a universal constant.
\end{lemma}
{\it Proof:} We note that the function  $v(x)-\frac{\alpha}{2n}|x-y_0|^2$ is super-harmonic in $B_r(y_0)$. The lemma follows from the mean value formula.
\qed
\vspace{10pt}

Now, we recall that $\vae(y_i)=0$ and $\Delta \vae \leq \alpha$ in $B_{1/4}(y_i)$.
So 
$$ \frac{1}{B_{1/4}}\int_{B_{1/4}(y_i)} \vae(x)\, dx \leq \vae (y_i) + \alpha C(n) $$
In particular, we have 
$$ u(y_i) \leq \int_{B_{1/4}(y_i)} \vae (x)\, dx \leq  C (\alpha,n)$$
\vspace{10pt}

\noindent {\it Step 4:} Steps 2 and 3 yield 
$$\sup_{B_i} u \leq C(\alpha,n)   (1+\alpha (\delta\eps^{-1})^2).$$
and since 
$ \Delta \vae \geq 0$ in $B_i\setminus \cap_{k\in\ZZ^n}\{k\}$, we have:
$$\vae (y) \leq  \frac{1}{B_{1/8}}\int_{B_{1/8}(y)} \vae(x)\, dx \leq C u(y)$$
for all $y\in B_i\setminus \cap_{k\in\ZZ^n} B_{1/4}(k)$.

\vspace{10pt}

It  follows that for every $\delta$ and for $\eps$ small enough, we have:
 $$
 \sup_{B_{\eps^{-1}} \setminus\cup_{k\in\ZZ^n} B_{1/4}(k)} \vae \leq C\delta^2 \eps^{-2}.
 $$
The definition of $\vae$ and the fact that $h_{\alpha,k} \geq 0 $ on $\pa B_{1/2}$  implies that
$$
\vae(x) \leq h_{\alpha,k}(x) + C\delta^2 \eps^{-2} \qquad \mbox{ in } B_{1/2}(k)
$$
for all $k\in\ZZ^n$. 
\qed
\vspace{30pt}

We now want to use the solution (\ref{eq:ovae}) of the obstacle problem (\ref{def:ov}) with $A=D$ to study the properties of  the free solution $w^\eps_0$ of
$$
\left\{
\begin{array}{ll} 
\Delta w^\eps_0 = \alpha_0 - \sum_{k\in \ZZ^{n}\cap D} \gamma(k,\omega) \delta(x-\eps k) , & \mbox{ in } D
\\
w^\eps_0 =0  & \mbox{ on } \pa D
 \end{array}
 \right.
$$
We prove:
\begin{lemma}\label{lem:woe}
For every $k\in\ZZ^n$, $w^\eps_0$ satisfies
\begin{equation}\label{eq:w0}
h_{\alpha,k}^\eps (x) -o(1)\; \leq\;  w^\eps_0 (x)\; \leq\; h_{\alpha,k}^\eps (x) +o(1)\quad  \forall x\in  B_{\eps/2}(\eps k)\cap \A\quad \mbox{ a.e. } \omega\in\Omega,
\end{equation}
with
$$ h_{\alpha,k}^\eps (x):= 
\left\{
\begin{array}{ll}
\ds
\frac{\alpha_0}{2n}|x-\eps k|^2+\frac{\eps^n r(k,\omega)^{n-2}}{|x-\eps k|^{n-2}}
 & \mbox{ if } n\geq 3 \\[10pt]
\ds\frac{\alpha_0}{2n}|x-\eps k|^2 -r(k,\omega) \eps^2 \log |x-\eps k|& \mbox{ if } n=2,
 \end{array}\right.
$$
In particular: 
\begin{equation}\label{eq:woe1}
w^\eps_0(x) = 1+o(1)\quad \mbox{ on } \pa T_\eps\cap \A
\end{equation}
\end{lemma}

Note that with this definition of $h_{\alpha,k}^\eps$, we have $h_{\alpha,k}^\eps (x)= \eps^2 h_{\alpha,k}(x/\eps) $ for $n\geq 3$ and $h_{\alpha,k}^\eps (x)= \eps^2 h_{\alpha,k}(x/\eps) +r\eps^2 \log \eps $ for $n=2$.
\vspace{10pt}

{\bf Proof.} For every $\alpha$, we denote by $\overline w^\eps_\alpha$  the function
$$ \overline w^\eps_\alpha (x)= \eps^2 \; \ov_{\alpha,\eps^{-1}\A}(x/\eps) ,$$
defined in $\A$ and satisfying $\overline w^\eps_\alpha = 0 $ on $\pa \A$.

\begin{enumerate}
\item For every $\alpha> \alpha_0$, we have
$$ \Delta (w^\eps_0-w^\eps_\alpha) \geq \alpha_0-\alpha 
$$
and $w^\eps_0-w^\eps_\alpha=0$ on $\pa \A$. This implies
$$ w^\eps_0(x_0)-w^\eps_\alpha(x_0) \leq \int_\A G(x_0,x) (\alpha_0-\alpha)\, dx$$
where $G(\cdot,\cdot)$ is the Green function on $\A$ ($\Delta G = \delta_{x_0}$ and $G=0$ on $\pa \A$).
Note that we have
$$G(x_0,x) \geq -h(x-x_0) \qquad \forall x,\, x_0 \in  \A ,$$
and so
$$
 w^\eps_0(x_0)-w^\eps_\alpha(x_0) \leq (\alpha - \alpha_0) \int_\A  h(x-x_0) \, dx.$$
We deduce
$$ \sup_\A  \left( w^\eps_0-w^\eps_\alpha \right) \leq \left\{\begin{array}{ll}
C|\A |^{1/(n-1)} \rho_D  \,|\alpha-\alpha_0| &\mbox{ if } n\geq 3 \\ 
C|\A | \rho_D \log \rho_D  \,|\alpha-\alpha_0|  &\mbox{ if } n=2 ,
\end{array}
\right.
$$
with
$$
\rho_D = \inf \{\rho\, ;\, D\subset B_\rho\}.
$$
Hence we have
$$ w^\eps_0\leq w^\eps_\alpha +O(\alpha-\alpha_0).$$
Using Lemma \ref{lem:1.2} (ii) (since $\alpha> \alpha_0$), we deduce:
$$ w^\eps_0 \leq h_{\alpha,k}^\eps (x) +O(\alpha-\alpha_0)+o(1)\quad  \forall x\in  B_{\eps/2}(\eps k)\quad \mbox{ a.e. } \omega\in\Omega .$$ 
which gives the second inequality in (\ref{eq:w0}).

\item Similarly, we observe that for every $\alpha \leq \alpha_0$, we have
$$ \Delta (w^\eps_\alpha-w^\eps_0) \geq \alpha-\alpha_0 -\alpha 1_{\{w^\eps_\alpha = 0\}} .
$$
Proceeding as before, we deduce that for $n\geq 3$, 
$$ \sup_\A  \left( w^\eps_\alpha-w^\eps_0  \right) \leq C\rho_D \left[|\A |^{1/(n-1)}(\alpha_0-\alpha) + C \alpha |\{w^\eps_\alpha = 0\}|^{1/(n-1)}\right]$$
and a similar inequality for $n=2$.
Using Lemma \ref{lem:1.2} (i), we get
$$ w^\eps_0 \geq h_{\alpha,k}^\eps - o(\eps^2) - O(\alpha_0-\alpha) - C \alpha |\{w^\eps_\alpha = 0\}|^{1/(n-1)}.$$
Finally, since
$$\lim_{\eps\rightarrow 0} |\{w^\eps_\alpha = 0\}| = 0 $$
for all $\alpha \leq \alpha_0$, and (\ref{eq:w0}) follows.
\end{enumerate}
\qed

\vspace{10pt}

\subsection{Proof of Proposition \ref{prop:1}} 

We are now in position to complete the proof of Proposition \ref{prop:1}:
We define
$$
w^\eps (x,\omega) = \inf\left\{w(x)\, ; \, \Delta w \leq \alpha_0 \mbox{ in } D\setminus T_\eps \, ,\,
\begin{array}{l} 
w \geq 1  \mbox{ on } T_\eps \cap D \\
w = 0  \mbox{ on } \pa D \setminus T_\eps 
 \end{array}
 \right\},
$$
it is readily seen that
$$ 
\left\{
\begin{array}{ll}
w^\eps (x,\omega)  = 1 & \mbox{ on } \pa T_\eps ,\\[5pt]
\Delta w^\eps (x,\omega)  = \alpha_0 & \mbox{ on }  D\setminus T_\eps,\\[5pt]
w^\eps (x,\omega)  = 0 & \mbox{ on } \pa D\setminus T_\eps .
\end{array}
\right.
$$
So in order to complete the proof, we only have to show that $w^\eps \longrightarrow 0$ in $H^1(D)$-weak as $\eps$ goes to zero.
More precisely, we will show that $w^\eps $ converges to zero in $L^p$ strong and is bounded in $H^1$.
\vspace{10pt}

\noindent {\bf Strong convergence in $L^p$:}\\
First of all, (\ref{eq:woe1}) yields
$$
w^\eps_{0}(x) -o(1)\; \leq\;  w^\eps  (x,\omega)\; \leq\; w^\eps_{0}(x) +o(1)\quad  \forall x\in  \De \quad \mbox{ a.e. } \omega\in\Omega,
$$
which in turns imply (using Lemma \ref{lem:woe}  again):
\begin{equation}\label{owh}
h_{\alpha,k}^\eps (x) -o(1)\; \leq\;  w^\eps  (x,\omega)\; \leq\; h_{\alpha,k}^\eps (x) +o(1)\quad  \forall x\in  B_{\eps/2}(\eps k)\quad \mbox{ a.e. } \omega\in\Omega.
\end{equation}
Next, a simple computation shows that
$$
\int_{B_\eps\setminus B_{a^\eps}} |h^\eps_{\alpha,k}|^p\, dx \leq \left\{
\begin{array}{ll}
C \eps^n\left( \eps^{\frac{2n}{n-2}} +\eps^{2p}\right) & \mbox{ if } n\geq 3,\\[8pt]
C\eps^2 \eps^{2p} (\log \eps)^p & \mbox{ if } n=2
\end{array}
\right.
$$
Since $\#\{\eps\ZZ^n \cap D\}\leq C\eps^n$ for all $n$, we deduce from (\ref{owh}) that
\begin{equation}\label{eq:Lp2} 
||w^\eps ||_{L^p} \leq \left\{
\begin{array}{ll}
C\left( \eps^{\frac{2n}{p(n-2)}}+\eps^2\right)  & \mbox{ if } n\geq 3, \\[5pt]
C\eps^2(\log\eps)\, & \mbox{ if } n= 2.
\end{array}
\right.
\end{equation}
In particular
$$
w^\eps \longrightarrow 0 \quad \mbox{  in } L^p-\mbox{strong,  for all $p\in[1,\infty)$}.
$$
\vspace{10pt}

\noindent {\bf Bound in $H^1$:}\\ 
First of all, a simple integration  by parts together with the fact that $w^\eps = 1$ on $\pa T_\eps$ yields
$$ \int_{ D_\eps} |\na w^\eps|^2 \, dx \leq \alpha_0 |D | + \int_{\pa T_\eps  } |\na w^\eps| d\sigma(x)$$
where $\pa T_\eps  = \cup_{} \pa S_\eps(k,\omega) $.
So we need an estimate in $\na w^\eps$ along $\pa S_\eps(k,\omega) = \pa B_{a^\eps(r(k,\omega))}$.

We consider the function
$$ z(x) = \left\{
\begin{array}{ll}
w^\eps  (x) - h_{\alpha,k}^\eps (x) +\frac{\alpha_0}{2n}r^2  \eps^{n/(n-2)} & \mbox{ when } n\geq 3 \\ [8pt]
w^\eps  (x) - h_{\alpha,k}^\eps (x) +\frac{\alpha_0}{2n} r^2 e^{2\frac{\eps^{-2}}{r}  }& \mbox{ when } n =2.
\end{array}
\right.
$$
It satisfies
$$
\left\{
\begin{array}{ll}
\Delta z=0 & \mbox{  in } B_{1/2} (\eps k)\setminus B_{ a^\eps(r(k,\omega))}(\eps k), \\ [5pt]
z(x)=o(1)\quad&  \mbox{ in } B_{1/2} (\eps k)\setminus B_{ a^\eps(r(k,\omega))}(\eps k) \\[5pt]
z(x)=0 & \mbox{ along } \pa B_{ a^\eps(r(k,\omega))}(\eps k) ,
\end{array}
\right.
$$
and so
$$ |\na z(x)| \leq 
\left\{
\begin{array}{ll}
o(r^{n-2} \eps^n\eps^{-\frac{n(n-1)}{n-2}}) = o(\eps^n a^\eps(r)^{-(n-1)}) & \mbox{ if } n\geq 3 ,\\
o( \eps^2e^{r^{-1}\eps^{-2}}) = o(\eps^n a^\eps(r)^{-(n-1)})
& \mbox{ if } n=2.
\end{array}
\right.
$$
on $\pa B_{ a^\eps(r(k,\omega))}(\eps k) $.
It follows that
$$ |\na w^\eps | \leq |\na h_{\alpha,k}^\eps (x)| +|\na z(x)| \leq C \eps^n a^\eps(r(k,\omega))^{-(n-1)}$$
along $\pa B_{ a^\eps(r(k,\omega))}(\eps k)  $,

We deduce
\begin{eqnarray*}
  \int_{ D_\eps} |\na w^\eps|^2 \, dx 
&\leq &  \alpha_0 |D | + \int_{\pa T_\eps  } |\na w^\eps| d\sigma(x)\\
&\leq &   \alpha_0 |D | + \sum_{k\in\ZZ^n\cap\eps^{-1}D} \int_{\pa B_{a^\eps(r(k,\omega))}(\eps k) } |\na w^\eps| d\sigma(x) \\
&\leq &  \alpha_0 |D | + C \eps^{-n} a^\eps(\overline r)^{n-1}  \eps^n a^\eps(\overline r)^{-(n-1)} \\
&\leq &   C,
\end{eqnarray*}
and the proof is complete.
\qed

\vspace{20pt}

\section{Proof of Proposition \ref{prop:1}: General case}\label{sec:general}
In this section, we treat the case where the  sets $S_\eps (k,\omega)$ have unspecified shape, but satisfy Assumption 2:
$$
\capa(S_\eps(k,\omega)) = \eps^n \gamma(k,\omega).
$$
Throughout this section we assume $n\geq 3$.

The proof makes use of the result of the previous section, after noticing that away from $\eps k$, the hole $S_\eps(k,\omega)$ is equivalent to a ball of radius $a^\eps(r(k,\omega))$, where
$$
a^\eps(r)=r \eps^{n/(n-2)}, \qquad r(k,\omega)=\left(\frac{\gamma(k,\omega)}{n(n-2)\omega_n}\right)^{1/(n-2)}
$$
More precisely, we will rely on the following lemma:
\begin{lemma}\label{lem:cap}
For any $k\in\ZZ^n$ and $\omega\in\Omega$, 
let  $\vphi^\eps_{k}(x,\omega)$ be defined by
$$  \vphi^\eps_{k}(x,\omega)= \inf \left\{v(x)\, ;\,  \Delta v \leq 0  \, ,\left\{ \!\! \begin{array}{l}v(x)\geq 1, \quad \forall x\in S_{\eps}(k,\omega)\\[5pt]
\lim_{|x|\rightarrow \infty} v(x) = 0 \end{array}\right.\right\}
$$
Then for any $\delta>0$, there exists $R_\delta$ such that
$$ \left| \vphi^\eps_{k}(x,\omega) - \eps^n \gamma (k,\omega) h(x-\eps k)\right| \leq \delta \eps^n  h(x-\eps k) \quad$$
for all $x$ such that $|x-\eps k| \geq a^\eps(R_\delta)$
and for all $\eps>0$.

Moreover, $R_\delta$ depends only on the constant $M$ appearing in Assumption~{\bf1}. In particular, $R_\delta$ is independent on $k$ and $\omega$.
\end{lemma}

\begin{enumerate}
\item For a given  $\delta>0$, Lemma \ref{lem:cap} implies that for every $k\in\ZZ^n$ and $\omega\in\Omega$ there exists a constant $R_\delta(k,\omega)$ such that
\begin{equation}\label{eq:gk}
\left| \vphi^\eps_k  (x,\omega) - \frac{\eps^n r(k,\omega)^{n-2}}{|x-\eps k|^{n-2}}  \right| \leq \delta \left(\frac{r}{R_\delta}\right)^{n-2}  \quad \mbox{ in }   B_{2a^\eps(R_\delta)}\setminus B_{a^\eps(R_\delta)} (\eps k)
\end{equation}
for all $\eps>0$.
Moreover,  it is readily seen that for any $R$ there exists $\eps_1(R )$ such that 
\begin{equation}\label{eq:Rd}
a^\eps(R) \leq \eps^\sigma /4\quad  \mbox{ for all $\eps\leq \eps_1$.}
\end{equation}
for some $\sigma >1$.
Finally, we note that by definition of $\vphi^\eps_k$, we have
\begin{equation} \label{eq:gradg}
\int_{\RR^n} |\na \vphi^\eps_k|^2\, dx = \capa (S_\eps(k)) =\eps^{n} \gamma(k,\omega) 
\end{equation}

\vspace{10pt}

\item Next, let $\alpha_0$  and $w^\eps $ be the coefficient and corresponding corrector constructed in the previous section, and associated with holes $S_\eps$ of radius $r(k,\omega)$.
Lemma \ref{lem:woe} implies that for $\delta$ and $R$ given, there exists $\eps_2(\delta,R)<\eps_1(R)$ such that
for all $\eps\leq \eps_2(\delta,R)$, we have
\begin{equation}\label{eq:we}
\left|w^\eps (x) - \frac{\eps^n r(k,\omega)^{n-2} }{|x-\eps k|^{n-2}} \right| \leq \frac{\delta}{R^{n-2}} \qquad \mbox{ in }  B_{\eps/2} (\eps k),
\end{equation}
in dimension $n\geq 3$.
Note that thanks to (\ref{eq:Rd}), Inequality (\ref{eq:we})  holds in particular in
$ B_{2a^\eps(R)}\setminus B_{a^\eps(R)} (\eps k)$.
\vspace{10pt}

\end{enumerate}

The  corrector  given by Proposition \ref{prop:1} will be constructed by gluing together the functions $\vphi^\eps_k$ (near the holes $S_\eps(k)$ and the function $w^\eps $ (away from the holes).
The gluing will have to be done in a very careful way so that the corrector satisfies all the properties listed in Proposition   \ref{prop:1}:
For a given $\eps$, we define $\delta_\eps$ to be the smallest positive number such that (\ref{eq:Rd}) and (\ref{eq:we}) hold with $\delta=\delta_\eps$ and $R=R_{\delta_\eps}$.
From the remarks above, we see that $\delta_\eps$ is well defined as soon as $\eps$ is small enough (say smaller than $\eps_2(1,R_1)$).
Moreover, 
for any $\delta>0$, there exists $\eps_0=\eps_2(\delta,R_\delta)$ such that
$$\delta_\eps \leq \delta \qquad \forall \eps\leq \eps_0.$$
In particular 
$$\lim_{\eps\rightarrow 0} \delta_\eps = 0.$$
From now on, we write
$$ R_\eps  = R_{\delta_\eps}.$$
\vspace{20pt}

We are now ready to define the corrector $\cor$ :
Let $\eta_\eps(x)$ be a function defined on $D$ such that
$$ \begin{array}{ll}
\eta_\eps(x) = 1 &\mbox{ on } D \setminus \cup_{k\in\ZZ^n} B_{2 a^\eps(R_\eps)}(\eps k)\\[5pt]
\eta_\eps(x)=0&\mbox{  on }  \cup_{k\in\ZZ^n}  B_{a^\eps(R_\eps)}(\eps k).
\end{array}
$$
and satisfying 
$$|\na \eta_\eps|\leq C  a^\eps(R_\eps)^{-1}\quad\mbox{  and }\quad |\Delta \eta_\eps|\leq Ca^\eps(R_\eps)^{-2} $$
in $B_{2 a^\eps(R_\eps)}\setminus  B_{a^\eps(R_\eps)}(k).$
We then define $\cor(x,\omega)$ in $D$ by:
$$\cor (x,\omega)  =   \eta_\eps (x) w^\eps (x,\omega )+(1-\eta_\eps(x))  \sum_{k\in \ZZ^n\cap D} \vphi^\eps_{k} (x,\omega)\,1_{B_{\eps/2}(\eps k)}(x) .$$
It satisfies
$$
\cor (x,\omega)=\left\{
\begin{array}{ll}
\vphi^\eps_{k} (x)& \mbox{ in }   B_{2a^\eps(R_\eps)}(k)\setminus S_\eps(k) \quad \forall k\in\ZZ^n \\[5pt]
w^\eps (x) &\ds  \mbox{ in } D \setminus \cup_{k\in\ZZ^n} B_{a^\eps(R_\eps)}.
\end{array}
\right.
$$
To simplify the notations in the sequel, we denote
$$ 
\vphi^\eps (x):= 
 \sum_{k\in \ZZ^n\cap D} \vphi^\eps_{k} (x,\omega)\,1_{B_{\eps/2}(\eps k)}(x)
$$

The properties of $\cor$ are summarize in the following lemma, which implies Proposition \ref{prop:1} with (\ref{H1'}) instead of the first equation:
\begin{lemma}
The function $\cor$ satisfies
\item [(i)] $ \cor = 1$ on $S_\eps$ for any $\eps>0$.
\item [(ii)] $ \cor$ converges to zero as $\eps$ goes to zero in $L^p(D)$ strong for all $p\in[2,\infty)$ and
$$ || \cor||_{L^p} \leq C \eps^\frac{2n }{p(n-2)}
\qquad \forall p\geq 2 $$
\item[(iii)] $ \cor$ is bounded in $H^1(D)$.
\item  [(iv)] $\Delta \cor$ converges to $\alpha_0$ in $L^1(D)$ and thus satisfies (\ref{H1'}).
\end{lemma}
{\it Proof:}
\begin{enumerate}
\item[(i)] Immediate consequence of the definition of $\cor$ since $\vphi^\eps_{k} = 1$ on $S_\eps(k,\omega)$.
\item[(ii)] 
Assumption $\bf 1$ yields
$$  \vphi^\eps_{k}(x,\omega) \leq C \eps^n \gamma (k,\omega) h(x-\eps k)$$
for all $x$ such that $|x-\eps k| \geq a^\eps(M)$.
Since $\vphi^\eps_{k}\leq 1$ in $B_{a^\eps(M)}(\eps k)$, we deduce: 
\begin{eqnarray*}
\left\| (1-\eta_\eps)  \vphi^\eps  \right\|^p_{L^p(\RR^n)} &  \leq &  \sum_{k\in \ZZ^n\cap \eps^{-1} D}  \left\|  \vphi^\eps_{k}\,1_{B_{\eps/2}(\eps k)}\right\|^p_{L^p({\cup B_{R(k)a(\eps)}})} \\
&  \leq &   \sum_{k\in \ZZ^n\cap\eps^{-1}  D} \int _{B_{a^\eps(M)(\eps k)}}\!\!\!\!\!\!\!\!\!\!\!\!  (\vphi^\eps_{k} (x))^p\, dx  \\
&& \quad + C \sum_{k\in \ZZ^n\cap\eps^{-1}  D} \int _{B_{2a^\eps(R_\eps)(\eps k)}}\!\!\!\!\!\!\!\!\!\!\!\!  ( \eps^n \gamma (k) h(x-\eps k))^p \, dx \\
&\leq  & \sum_{k\in \ZZ^n\cap\eps^{-1}  D} a^\eps(M)^n  \\
&& \quad + C \overline \gamma\sum_{k\in \ZZ^n\cap\eps^{-1}  D}  \eps^{pn}  (a^\eps(R_\eps))^{n-p(n-2)} 
\end{eqnarray*}
Using (\ref{eq:Rd}) and the definition of $a(\eps)$, we deduce:
\begin{eqnarray*}
\left\| (1-\eta_\eps) \vphi^\eps \right\|^p_{L^p(\RR^n)} 
&\leq & C   \eps^{-n}M^n\eps^\frac{n^2}{n-2} + C \overline \gamma\sum_{k\in \ZZ^n\cap D}  \eps^{pn}  \eps^{n-p(n-2)} 
\\
&\leq &C   M^n \eps^\frac{2n }{n-2} + C \overline \gamma\sum_{k\in \ZZ^n\cap D}  \eps^{n+2p}
\\
& \leq &C   M^n \eps^\frac{2n }{n-2} + C \overline \gamma  \eps^{2p}
\end{eqnarray*}
where $2p\geq \frac{2n }{n-2} $ if $p\geq 2 $ and $n\geq 3$.

Using (\ref{eq:Lp2}), it follows that
\begin{eqnarray*}
|| \cor ||_{L^p(D)}  & \leq &  || w^\eps  ||_{L^p(D)} + C\left(  \eps^\frac{2n }{n-2} \right)^{1/p}   \\
&\leq & C \eps^\frac{2n }{p(n-2)}
\end{eqnarray*}
for all  $p\geq 2$.

\item[(iii)] Next, we want to show that $\cor$ is bounded in $H^1(\De)$. 
First, we note that in $B_{\eps/2}(\eps k)$, we have:
\begin{equation}\label{eq:grad}
\nabla \cor = \na \eta_\eps (w^\eps -\vphi^\eps_{k}) +  \eta_\eps \na\overline  w^\eps_0+(1-\eta_\eps) \na \vphi^\eps_{k} 
\end{equation}
where the function $\na \eta_\eps$ is supported in $B_{2a^\eps(R_\eps)}(\eps k)\setminus B_{a^\eps(R_\eps)}(\eps k) $ and satisfies  
$$|\na \eta_\eps|\leq C (a^\eps(R))^{-1}.$$
Since 
$|w^\eps -\vphi^\eps_{k} | \leq C\frac{\delta_\eps}{R_\eps^{n-2}}$  in $B_{2a^\eps(R_\eps)}(\eps k)\setminus B_{a^\eps(R_\eps)}(\eps k) $, we deduce
\begin{eqnarray*}
\int_D | \na \eta_\eps (w^\eps -\vphi^\eps) |^2 \, dx & \leq & 
\sum_{k \in \eps\ZZ^n\cap D} \int_{B_{2 a^\eps(R_\eps)}(\eps k)} | \na \eta_\eps (w^\eps -\vphi^\eps_{k}) |^2 \, dx \\
& \leq & 
\sum_{k \in \eps\ZZ^n\cap D} ( a^\eps(R_\eps))^n(a^\eps(R_\eps))^{-2} \frac{\delta_\eps ^2}{R_\eps ^{2(n-2)}} \\
& \leq & 
\sum_{k \in \eps\ZZ^n\cap D} R_\eps ^{-(n-2)} \eps^n \delta_\eps ^2 \\
& \leq &C\eps^{-n} \eps^n = C,
\end{eqnarray*}
since we can always assume that $\delta_\eps<1$ and $R_\eps\geq 1$.
Finally, since $w^\eps $ and $\vphi^\eps $ are both bounded in $H^1$ (thanks to (\ref{eq:gradg})), (\ref{eq:grad}) implies
$$||\nabla \cor||_{L^2} \leq C.$$

\item[(iv)] It remains to evaluate the Laplacian of $\cor$. We have:
$$ 
\Delta \cor = \alpha  -(1- \eta_\eps) \alpha +2 \na \eta_\eps \cdot \na ( w^\eps -\vphi^\eps) + \Delta \eta_\eps\, ( w^\eps  - \vphi^\eps) \mbox{ in } \De.
$$
Moreover, (\ref{eq:we}) and (\ref{eq:gk}) yield
$$|w^\eps  - \vphi^\eps_{k} | \leq \frac{\delta_\eps }{R_\eps ^{n-2}}\qquad  \mbox{ in }B_{2a^\eps(R_\eps)}\setminus B_{a^\eps(R_\eps)} ,$$
and by definition of $w^\eps $ and $ \vphi^\eps_{k}$, we have
$$
\Delta( w^\eps  - \vphi^\eps_{k} -\frac{\alpha_0}{2n}|x-\eps k|^2) = 0 \qquad \mbox{ in } B_{4a^\eps(R)}\setminus B_{a^\eps(R_\eps)/2}.
$$
Interior gradient estimates thus implies
$$  | \na(w^\eps  - \vphi^\eps_{k})| \leq \frac{\delta_\eps }{R_\eps ^{n-2}} a^\eps(R_\eps)^{-1}  + C a^\eps(R_\eps)$$
in $B_{2a^\eps(R_\eps)} \setminus B_{a^\eps(R_\eps)} $.
We deduce (using (\ref{eq:Rd})):
\begin{eqnarray*}
&&\!\!\!\!\!\!\!\! \int_\De |\Delta \cor - \alpha|\, dx \\
& &\quad\quad \leq  \int_\De (1- \eta_\eps) \alpha \, dx +
\int_\De |\na \eta_\eps | |\na ( w^\eps -\vphi^\eps)| \, dx\\
&& \quad\qquad  + \int_\De |\Delta \eta_\eps| | w^\eps  - \vphi^\eps|\, dx \\
&& \quad\quad \leq  \sum_{k \in \eps\ZZ^n\cap \eps^{-1} D} a^\eps(R_\eps)^{n} \\
&& \quad\qquad+
\sum_{k \in \eps\ZZ^n\cap \eps^{-1} D}a^\eps(R_\eps)^{-1}\int_{ B_{2a^\eps(R_\eps)}\setminus B_{a^\eps(R)} } |\na ( w^\eps -\vphi^\eps_{k})| \, dx \\
&&\quad\quad \quad + \sum _{k \in \eps\ZZ^n\cap \eps^{-1}D}a^\eps(R_\eps)^{-2} \int_{ B_{2a^\eps(R_\eps)} \setminus B_{a^\eps(R)} }  | w^\eps_0 - \vphi^\eps_{k}|\, dx \\
& & \quad\quad\leq C  \sum_{k \in \eps\ZZ^n\cap \eps^{-1} D} a^\eps(R_\eps)^{n}+ C \sum _{k \in \eps\ZZ^n\cap\eps^{-1} D} 
\frac{\delta_\eps}{R^{n-2}} a^\eps(R_\eps)^{-2} (a^\eps(R))^{n} \\
& & \quad\quad\leq  C \eps^{-n} a^\eps(R_\eps)^{n}+ C  \delta_\eps \sum _{k \in \eps\ZZ^n\cap \eps^{-1}D} 
 \left(\frac{a^\eps(R_\eps)}{R_\eps}\right)^{n-2 } \\
& & \quad\quad\leq  C\eps^{(\sigma-1)n}+C \delta_\eps .
\end{eqnarray*}
In particular,
$$\lim_{\eps\rightarrow 0} \int_{\De} |\Delta \cor - \alpha_0 |\, dx \leq C\delta$$

\qed
\end{enumerate}

\appendix

\section{Proof of Lemma \ref{lem:cap}}
We recall that $n\geq 3$ in this section.
For any $k\in\ZZ^n$, we define
$\overline S_\eps(k) = \eps^{-\frac{n}{n-2}} S_\eps(k)$. Then Assumption {\bf 2} yields:
$$\capa (\overline S_\eps(k)) = \gamma(k)\leq\overline  \gamma.$$
and Assumption {\bf 1} gives
\begin{equation}\label{eq:M}
\overline S_\eps(k) \subset B_{M}(k).
\end{equation}
For the sake of simplicity, we take $k=0$.
We recall that $h$ is defined by 
$$
h(x)= \frac{1}{n(n-2)\omega_{n}}\frac{1}{|x|^{n-2}} .$$

Lemma  \ref{lem:cap} will be a consequence of the following lemma:
\begin{lemma} \label{lem:cap1}
Let  $\vphi$ be defined by
$$  \vphi(x)= \inf \left\{v(x)\, ;\,  \Delta v \leq 0  \, ,\left\{ \!\! \begin{array}{l} v(x)\geq 1, \quad \forall x\in \overline S_{\eps}(k,\omega) \\[5pt]
\lim_{|x|\rightarrow \infty} v(x) = 0 
\end{array}\right.\right\}
$$
Then for any $\delta>0$, there exists $R$, depending only on $\delta$ and $M$ such that
$$ \left| \vphi(x,\omega) -  \gamma  h(x)\right| \leq \delta h(x) \quad$$
for all $x$ such that $|x| \geq R$.
\end{lemma}

{\bf Proof:}
We recall that $\vphi$ solves
$$
\left\{
\begin{array}{ll}
\Delta \vphi(x) = 0\qquad & \mbox{ for all } x\in \RR^n\setminus S \\
\vphi(x) = 1 &  \mbox{ for all } x\in S\\
\lim_{|x| \rightarrow \infty } \vphi(x) = 0 .
\end{array}
\right.
$$
In particular, (\ref{eq:M}) and the maximum principle imply
\begin{equation}\label{eq:vphi} \vphi(x) \leq  M^{n-2} n(n-2)\omega_n h(x)=\frac{M^{n-2} }{|x|^{n-2}} \qquad\mbox{ in } \RR^n\setminus B_M(0).
\end{equation}

Next, we observe that
\begin{eqnarray*}
0  = -\int _{\RR^n\setminus S}\vphi\,  \Delta \vphi \, dx 
= \int_{\RR^n\setminus S} |\na \vphi|^2\,dx - \int_{\pa S}  \vphi\, \vphi_\nu\, d\sigma(x)  
\end{eqnarray*}
and so 
$$ \int_{\RR^n\setminus S} |\na \vphi|^2\,dx = \int_{\pa S}  \vphi\, \vphi_\nu\, d\sigma(x) = \int_{\pa S}  \, \vphi_\nu\, d\sigma(x)  .
$$
Moreover, for any $R\geq M$, we have
\begin{eqnarray*}
0  = \int _{B_R \setminus S}\,  \Delta \vphi \, dx 
=  \int_{\pa S}  \, \vphi_\nu\, d\sigma(x) +    \int_{\pa B_R}  \, \vphi_\nu\, d\sigma(x).
\end{eqnarray*}
We deduce:
\begin{equation}\label{eq:vphicapa} 
 \gamma = \int_{\RR^n\setminus S} |\na \vphi|^2\,dx = - \int_{\pa B_R}  \, \vphi_\nu\, d\sigma(x)  \quad  \mbox{ for all } R\geq M.
\end{equation}
\vspace{10pt}

We now introduce the function
$$\Theta(x) =h\left(\frac{x}{|x|^2}\right)^{-1} \vphi\left(\frac{x}{|x|^2}\right)=
n(n-2)\omega_{n} \frac{1}{|x|^{n-2}} \vphi\left(\frac{x}{|x|^2}\right)$$
defined for $x\in B_{1/M}(0)$.
A straightforward computation yields
$$\Delta \Theta = 0 \quad \mbox{ in } B_{1/M}(0)$$
and (\ref{eq:vphi}) implies
$$\Theta(x) \leq M^{n-2}n(n-2)\omega_n\quad \mbox{ in } B_{1/M}(0).$$
A more delicate computation, making use of the mean formula for harmonic functions, 
gives
 $$ 
\int_{\pa B_R}  \, \vphi_\nu\, d\sigma(x)  = -\Theta(0) .$$
Hence (\ref{eq:vphicapa}) yields
$$ 
\Theta(0) = \capa(\overline S_\eps) = \gamma$$

To conclude, we note that interior gradient estimates for harmonic functions imply the existence of a universal $C$ (depending only on $M$) such that
$$
|\Theta(x)-\gamma |\leq C|x| \qquad \mbox{ for all } |x|\leq 1/(2M).
$$
Inverting back, we deduce
$$
|\vphi(x) - \gamma h(x)| \leq \frac{C}{|x|}h(x) \mbox{ for all } |x|\geq 2M,
$$
which yields the result.
\qed
\bibliography{bibli}

\end{document}